# Preprint

Authors preprint version of a chapter to be published in the book *How theories work in mathematics education. Foundation, form and function*. H.L. Johnson, A. Bikner-Ahsbahs. Springer.

# Theory of Situations

Marie-Jeanne Perrin-Glorian, Laboratoire de Didactique André Revuz, Université Paris-cité, France, Orcid 0000-0001-9039-6528

Claire Margolinas, Laboratoire ACTé, Université Clermont-Auvergne, France, Orcid 0000-0001-5729-7966

## Abstract

Guy Brousseau began developing the theory of situations in the 1960s, with the ambition of building an experimental epistemology of mathematics. Over the decades, the theory has been structured, consolidated and extended, drawing on the contributions of numerous researchers. It is now recognized as one of the “genuine ground theories” in mathematics education. In this chapter, we will describe the foundations and reasons for the existence of this theory, its systemic approach to teaching and learning issues, its current architecture, which articulates a theory of mathematical situations and a theory of didactic situations, and its key concepts, as well as its functions, particularly its predictive and explanatory functions. We will illustrate this presentation with an example in geometry that also highlights the role played by the didactic engineering methodology in this theory and related work.

## Key words

## 1. Introduction

Niss and Jankvist (2022, p. 17) define a theory as “an organised network of concepts (including, ideas, notions, distinctions and terms) and claims (oftentimes taken as answers to questions) about some extensive domain of objects and phenomena” and they “grant the label ‘theory’ to [...] Brousseau’s (1997) theory of didactical situations.” This theory gradually developed and has been rooted in a scientific approach from the very beginning.

> *“Didactics, when it was set up as a field of research, was not constituted with normative ambitions, especially in France. [. . .] Research was built on the conviction that priority should be given to understanding the functioning of didactic systems, to clarify the processes of teaching and learning of mathematics which alone could be the basis of reasoned action.” (Artigue, in Gascòn & Nicolás, 2017, p. 4)*

The Theory of Situations (TS) was created by Guy Brousseau (1933–2024). During the 1960s and 1970s, he established its foundations with the project to study mathematics teaching scientifically, creating a new discipline he initially named “experimental epistemology of mathematics" (Brousseau, 2006, p. 10). To accomplish this, he needed to develop a methodological approach:

> *"I had to imagine a design that would put the activity of a class into relationship with a system of research that would satisfy a large number of conditions of very different orders. The three elements that seemed to me indispensable in order to make the act of teaching an object of scientific study are: observation, didactical engineering and the didactical study of mathematics itself." (Brousseau, 2006, p. 6).*

In 1973, Guy Brousseau's major methodological project, the COREM (Centre for the Observation and Research of Mathematical Teaching[1]), was established. It was a regular school that included resources for observation and collaboration with teachers for research purposes. During the 1980s and 1990s, these foundations gradually became consolidated and reorganized into a genuine theory, while maintaining the same main objectives and fundamental hypotheses. Numerous researchers around the world contributed to and continue to contribute to the development of the theory of situations and its implications. Many of them participated in person or gave video testimonials at the conference organized in tribute to Brousseau's work in July 2025 in Bordeaux (Derouet et al., 2025).

In this chapter, we present the foundations, main concepts, and functions of this theory. We begin with a brief overview of the TS and its foundation. We then move on to the various forms, structures and concepts that lie at the core of this theory. To illustrate some of its functions, we have chosen to develop the case of geometry, with two examples of didactic engineering: first the transition from practical to theoretical geometry at the transition from elementary school to middle school, as an example of a sequence of action situations, second the congruence of triangles as an example of combination of formulation situation and validation situation for the middle school.

# 2. Foundations

First, it seems important to contextualize the TS in relation to teaching and learning theories from the 1960s and 1970s. Then, we will distinguish between institutional knowledge and situational knowledge, which is an important part of the Theory of Mathematical Situations (TMS). Next, we will introduce the second part of the TS, the Theory of Didactic Situations in Mathematics (TDSM), and its main methodology, didactic engineering.

## 2.1. Learning as a combination of adaptation and acculturation

Brousseau (2000a) argues that psychologists have long demonstrated the importance of the environment for learning, which means that teaching cannot be viewed as simply the transmission of knowledge. He considers works based on the well-known theories of Piaget and Vygotsky very important as they draw attention to the role of personal actions and socio-cultural contexts for learning. However, he argues that teaching mathematics also involves an intentional process (teaching) within a specific institution (educational system) and a particular body of knowledge (mathematics) and thus for the necessity to consider the whole system of their interactions. A situation is a very broad concept characterized by a milieu and a stake. The milieu is the set of objects with which a subject interacts to achieve a stake. In some situations, prior knowledge may enable the subject to achieve the stake, but in others, new knowledge is required. In some situations, achieving the stake requires an exchange between two subjects or a debate among many. Brousseau emphasizes that in his project, the aim is to study situations not students. The concept of situation may be specified when we speak of a mathematical situation (see 2.3) or a didactic situation (see 2.4).

[1] See https://ardm.eu/guy-brousseau/guy-brousseau.com/le-corem/index.html and https://archivos-imac.uji.es/CRDM/index.php

For Brousseau, teaching coordinates two processes: adaptation to a milieu, coined in order to attain new knowledge, and acculturation to cultural knowledge (Fig. 1).

## 2.2. Situational knowledge and institutional knowledge

In French, Brousseau has drawn upon language in order to introduce a very important distinction between "connaissance" and "savoir": in English, Margolinas (2022) has proposed to translate those distinctions by "situational knowledge" and "institutional knowledge" (Margolinas & Bessot, 2024).

When trying to deal with a specific situation, someone may use or create situational knowledge to obtain the desired outcome. For an observer, this situational knowledge can sometimes be considered as mathematical knowledge. However, when a person is trying to act, this situational knowledge may be impossible to express verbally and may not be recognized as useful. In other situations (see formulation situation, validation situation in 3.2), the necessity of formulation and validation can gradually appear (Margolinas & Bessot, 2024).

In contrast, within the institution of mathematics, institutional knowledge is formalized, defined, and written in such a way as to be accepted by the community of mathematicians.

*Figure 1. The didactic situation in TS, schema inspired by Brousseau (2000a)*

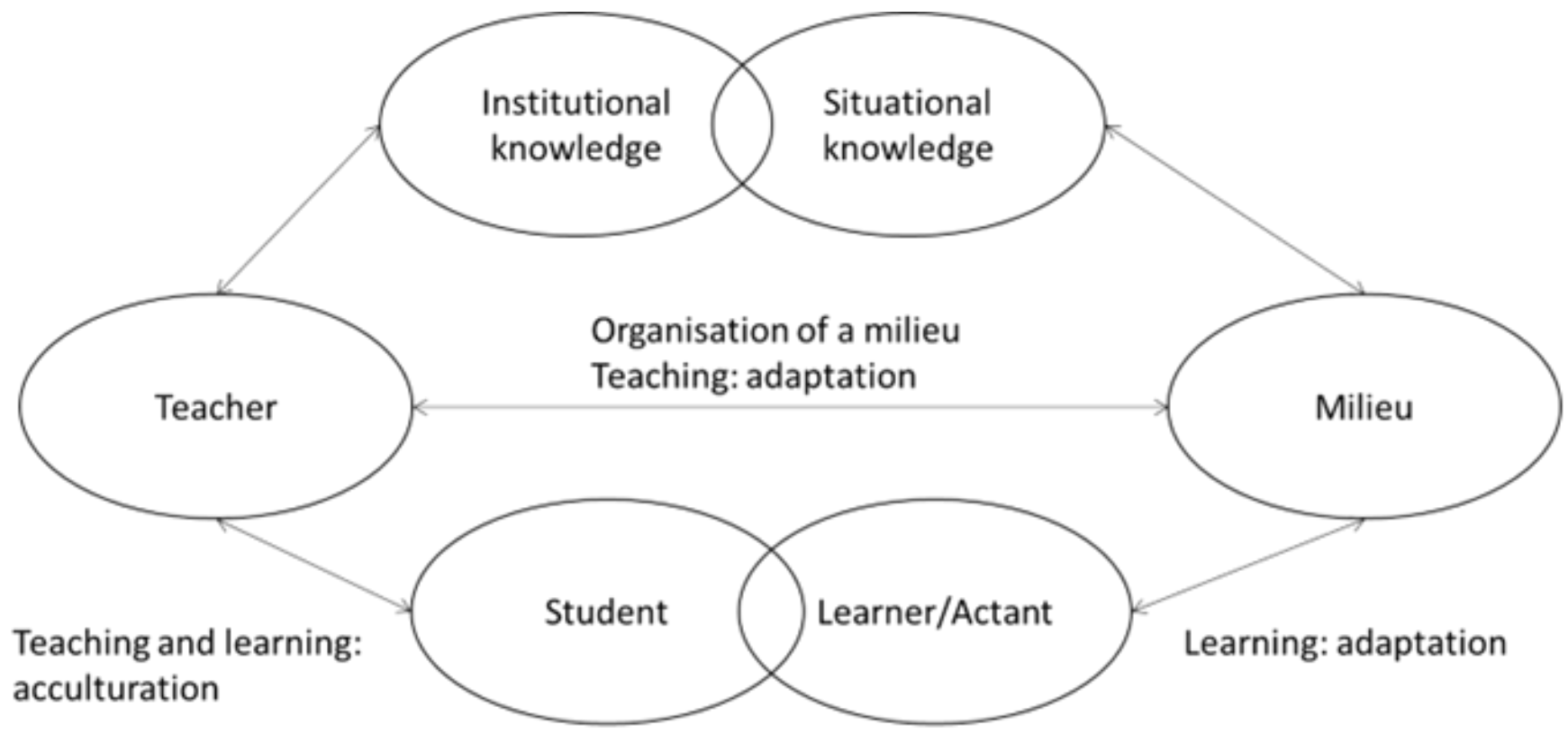

In a didactic situation, situational knowledge is used or created by a student acting in a situation. However, this situational knowledge must relate to institutional knowledge that the student must learn (Fig.1). This relationship is established partly by the organization of the situation and partly by the teacher's actions within the acculturation process (Margolinas & Bessot, 2024).

## 2.3. Theory of Mathematics Situations

In order to create situations and thus organize milieus and define stakes where mathematical knowledge might be necessary, Brousseau has been faced with the necessity of:

> *"redefining mathematics on the basis of the conditions of its appearance and use in human interactions and not directly. Every theorem needed to be associated with at least one situation that determined it as a correct solution but also offered incorrect options. These situations* *constitute models of the functioning of mathematics on which one could work a priori and predict their possibilities before realizing them. […] Inversely one could model non-experimental scholastic episodes observed in the course of ordinary classes and try to explain their logic." (Brousseau, 2006, p. 6).*

In this sense, the TS can be considered an epistemological revolution that views mathematics not only as a cultural given, but as a field of knowledge useful to society, whose usefulness must be rediscovered by describing it as knowledge that is useful in situations. This epistemological part of the TS has been called the TMS since the 2000's, although its dynamic was included in the TS right from the beginning.

At the ICME (International Congress on Mathematical Education) conference in Copenhagen, Brousseau (2004), as the Felix Klein medalist, introduced the principles of the TMS.

> *"A Theory of Mathematical Situations—Why?*
> *I made the assumption that*
> *– to every piece of mathematical knowledge there corresponds a collection of situations which can be resolved using this knowledge and reciprocally that*
> *– in any real environment of a student it is possible to choose elements of one or more situations that make it possible to identify the knowledge being brought into action by the student's actions." (Brousseau, 2004, p. 250)*

Guy Brousseau has devoted a great part of his work to reconsider mathematical knowledge in order to identify "the collection of situations which can be resolved using this knowledge" (Brousseau, 2004, p. 250), which he proposed to call a "fundamental situation." A fundamental situation (Bessot & Margolinas, 2025) must allow for the representation of all situations that represent the various meanings of a specific piece of knowledge using various didactic variables (see also section 5). This reconstruction of knowledge in the TMS, based on fundamental situations, is more than just a translation of mathematical knowledge in terms of situations; it also implies a new understanding of mathematical knowledge.

However, the TS inherently emphasizes observation and experience in real teaching contexts, as well as modeling interactions between teachers, students, and knowledge. Thus, it comprises another component presented in the next section.

## 2.4. Didactic engineering and Theory of Didactic Situations in Mathematics

The two components of the TS are intrinsically linked: the epistemological part, the TMS; and the study of learning and teaching mathematics, the Theory of Didactic Situations in Mathematics (TDSM), which includes a specific methodology: didactic engineering (Artigue, 2020; Margolinas & Bessot, 2024), see sections 4 and 5. In Guy Brousseau's work, these two components are always present. They appeared at the same time, from the beginning of the theory's creation even though the distinction was not made explicit until later.

The scope of the TDSM involves studying the relationship between knowledge, teachers and students (Brousseau et al., 2014). Close observation of these relationships in the classroom is required, which is why the COREM played an extremely important role in developing the TDSM. For some thirty years, Guy Brousseau's work included weekly work at the COREM: construction of mathematical situations with the COREM partners (teachers, researchers, etc.); observation of lessons followed by a discussion; data analysis (videos, statistics, etc.); not to mention daily discussions with his wife, Nadine Brousseau, who was a teacher at the school where the COREM was located. Even in an experimental setting, the intellectual conception of situations (the TMS component) needed to be confronted with the actual possibilities of teaching and learning (the TDSM component). Subsequently, other researchers (including the authors of this chapter) working within the TS framework focused their research on observing and analyzing teaching and learning in “regular” (non-experimental) classrooms (Hersant & Perrin-Glorian, 2005; Margolinas et al., 2005).

The originality of the TS consisted in not focusing on the individuals—teachers and students—or on the elementarization of mathematics, but rather in considering the system formed by teachers, students, and knowledge as they engage in a project defined by society. To conduct a scientific study of this system, it was necessary to develop a specific methodology: didactic engineering. The TS is therefore a theory that encompasses both fundamental and experimental research. Its purpose is not to improve teachers' or students' performance in standardized tests, but to understand how mathematics and the concepts behind it can be learnt and taught. We now come to the organization of the theory.

# 3. Form of the theory

The core concepts of the TS were elaborated and constantly reconstructed from the 1970s to the 1990s (Perrin-Glorian, 1994). We present here the present form, first the relationships between adidactic situation with its milieu and stake, didactic situations and didactic contract, second the main structures of adidactic situations and finally the two processes of devolution and institutionalization which are the essential ways for the teacher to make the adidactic situation exist inside the didactic situation.

## 3.1. Adidactic situation, didactic situation, milieu and didactic contract

In a didactic situation, the teacher has a didactic intention: to teach a specific mathematical concept to the students. Generally speaking, the usefulness of knowledge is revealed in situations that are not designed for teaching purposes, but that exist in everyday life, at work, or in mathematics, etc. The resulting didactic problem is then to create conditions that make this knowledge useful to students, considering their prior knowledge. The mathematical situation resulting from this research process in the TMS represents the most appropriate subject-milieu interactions. Brousseau (1997) named this mathematical situation an “*adidactic* situation[2]” when it is included in a didactic situation. It is the core of this didactic situation. The adidactic situation can be modeled as a game. In this game, the “actant” interacting with a milieu must try to achieve a defined stake using any available means within the established rules of the game (see the example on part 5, and more examples in Margolinas & Bessot, 2024). The term “actant” was introduced to refer to the actor in the game, whoever that may be. In the case of a successful devolution (see 3.4.2.), students can temporarily assume the role of an actant in relation to a milieu during their interactions with it.

In the early 80, Brousseau has also identified the didactic contract to model the mutual expectations of the teacher and students regarding knowledge: “The interpretation of the questions asked, the information provided, and the constraints imposed [...] are constants in the teacher's way of teaching. These (specific) habits of the teacher expected by the student and the behavior of the students expected by the teacher, this is the didactic contract" (Brousseau, 1980, translated by Margolinas & Bessot, 2024). The didactic contract is the result of a "negotiation" that is largely implicit. Figure 2 schematizes the structure of the didactic situation. In the TDSM, the adidactic situation is the core of the research for meaningful mathematics teaching.

[2] Brousseau originally spelled it “a-didactic,” but we use the spelling that subsequently became prevalent: adidactic.

*Figure 2. A model of the didactic situation*

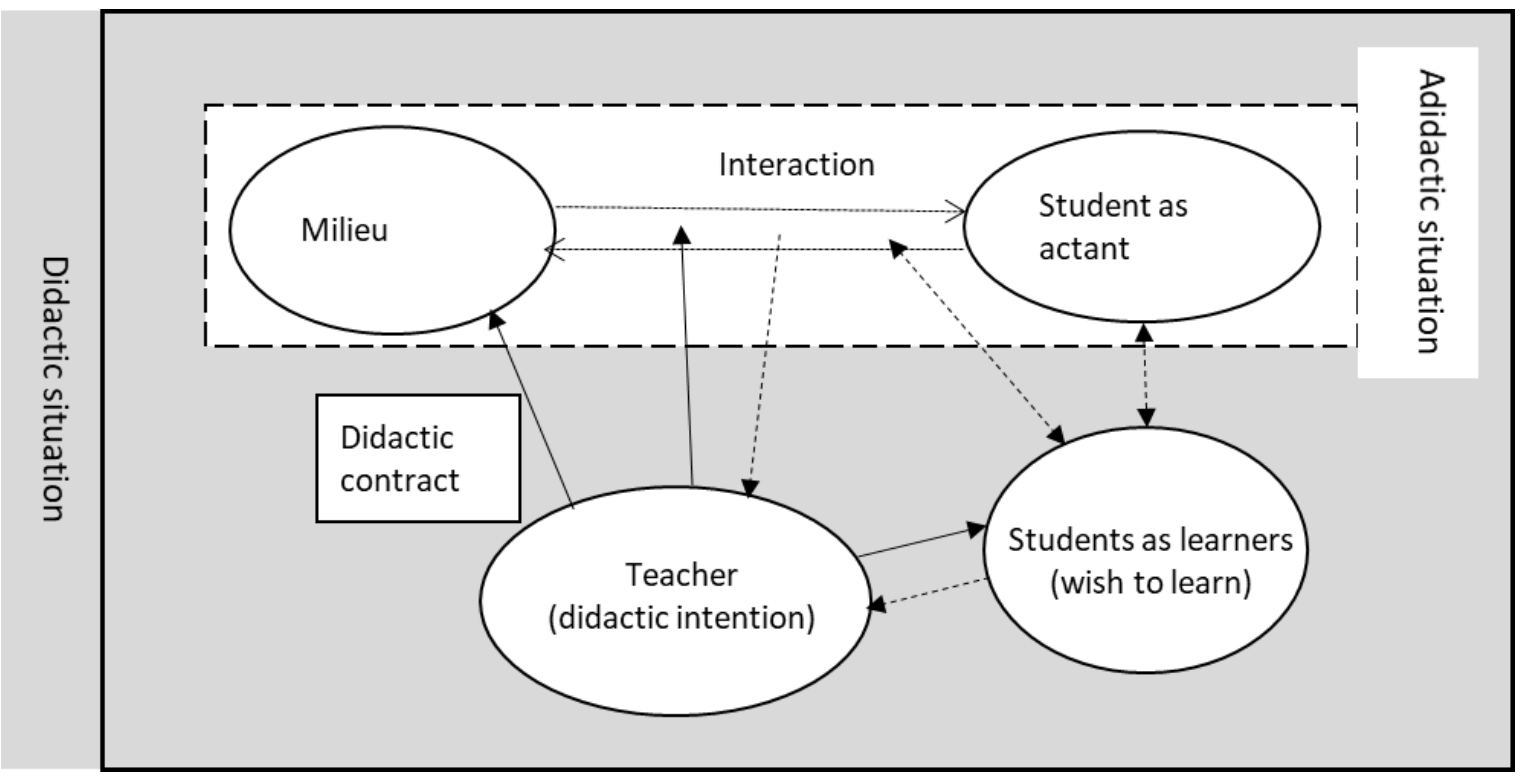

Brousseau stressed the paradoxical nature of the didactic contract and described several paradoxes. The various paradoxes (Brousseau, 1997, p. 40-45) originate from the first one, which is as follows: the teacher must obtain answers and actions from students that will lead to knowledge, but if the teacher says what they want the students to do, the students can only carry it out as an order, and not through the application of their own knowledge and judgment.

By the end of the 1980s, Brousseau began formalizing the structuring of the milieu to account for the "depth" of situations. This model was subsequently completed and modified (Bloch, 2002; Margolinas, 2002; Margolinas & Bessot, 2025; Salin, 2002). Due to space limitations, we cannot present this structure in detail here. Readers are encouraged to consult the references for more information. We can say, at least, that it considers the teacher's point of view, from the knowledge to be taught to the observation of students' procedures, as well as the student's point of view, from basic interaction with the milieu to reflection on acquired knowledge.

## 3.2. Different structures of adidactic situations

The distinction between situations of action, formulation, and validation was made by Brousseau as early as 1970 (Brousseau, 1997, p. 6-17) and refers to three functions of mathematical knowledge: knowledge needed to make a decision, to formulate, and to validate. To illustrate the definitions of these situations, we outline the key elements that we will develop in more depth in Section 5.

### 3.2.1. Action Situation

In an action situation (Fig. 3), someone must interact with the milieu to reach a determined stake. For example, say students are asked to draw a figure (action), from a starting figure using some tools. The milieu consists of the sheet of paper which includes the starting figure and the instruments. The drawn figure must correspond to a template traced onto a transparent sheet (stake). The students use or create situational knowledge (Fig. 2 above). See section 5.1 for a more detailed presentation of this example.

*Figure 3. Schema of action situation*

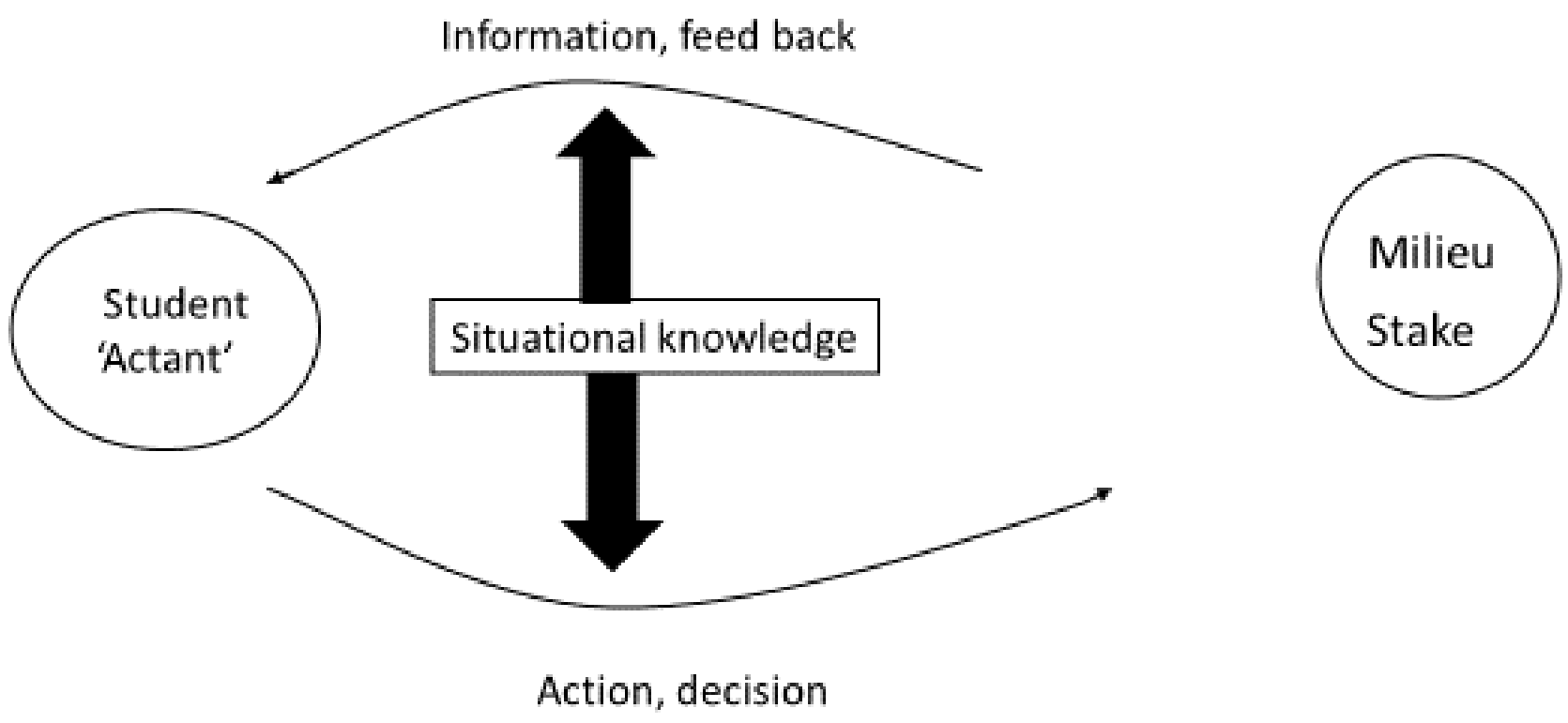

### 3.2.2. Formulation Situation

In a formulation situation (Fig. 4), the action on the milieu depends on information that will be provided later or by someone else. For example, it may involve two people with different access to actions and information (first type), or only one person (second type) who has access to information on one day but can only act on another day (e.g. a shopping list). In section 5.2, we present a formulation situation of the first type. Students are given information about a triangle and must write a message that allows another student to construct the same triangle (formulation situation). Finally, they compare the two triangles (stake).

*Figure 4. Schema of formulation situation*

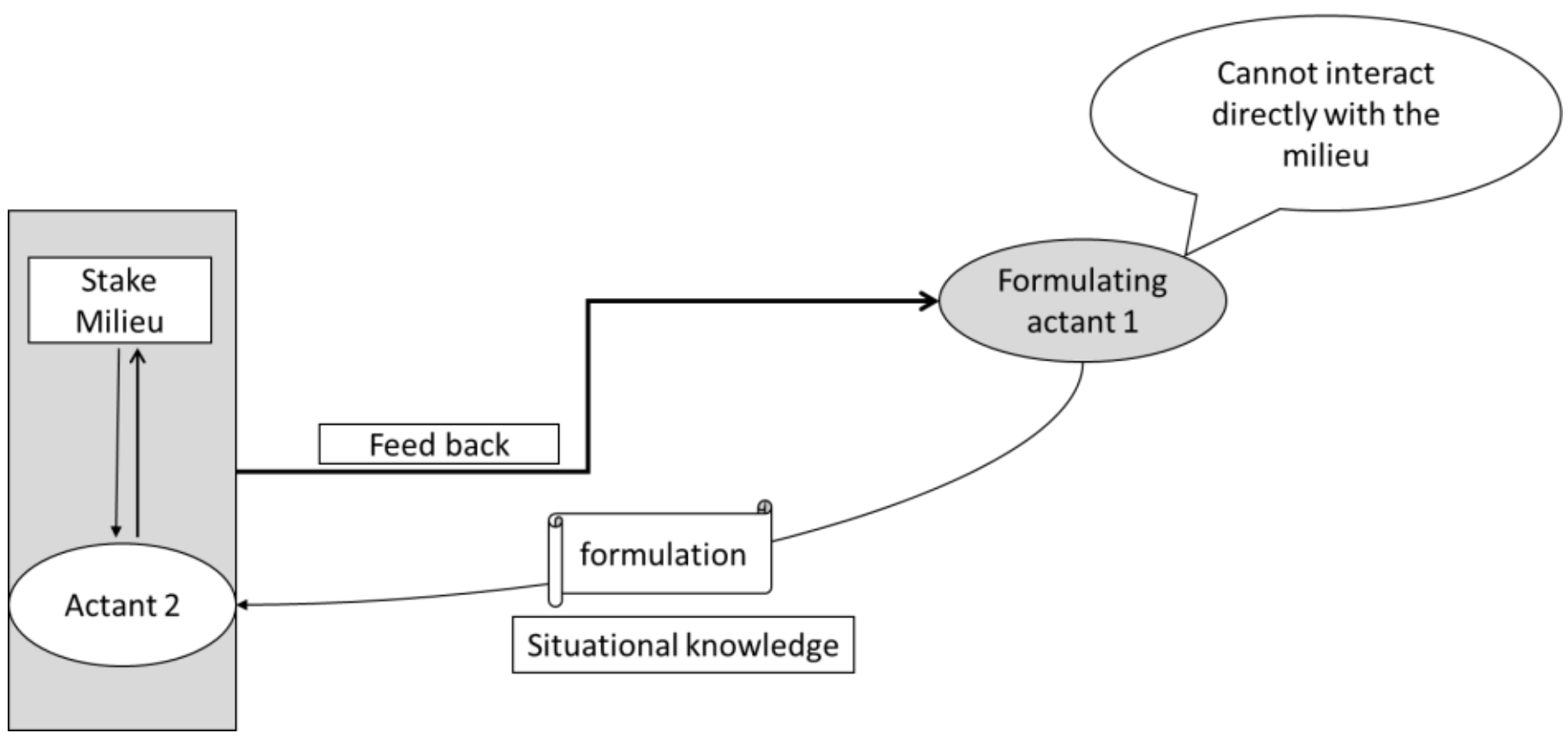

### 3.2.3. Validation Situation

A validation situation (Fig. 5) is a situation that seeks to create conditions for mathematical proof reasoning, Similarly, it includes an action situation and a formulation situation. Previously produced formulations will be useful for generating assertions whose validity will be questioned and tested through action. A proponent is an actant who formulates a statement that they consider true. This statement is submitted to the approval or rejection of the opponents.

*Figure 5. Schema of validation situation*

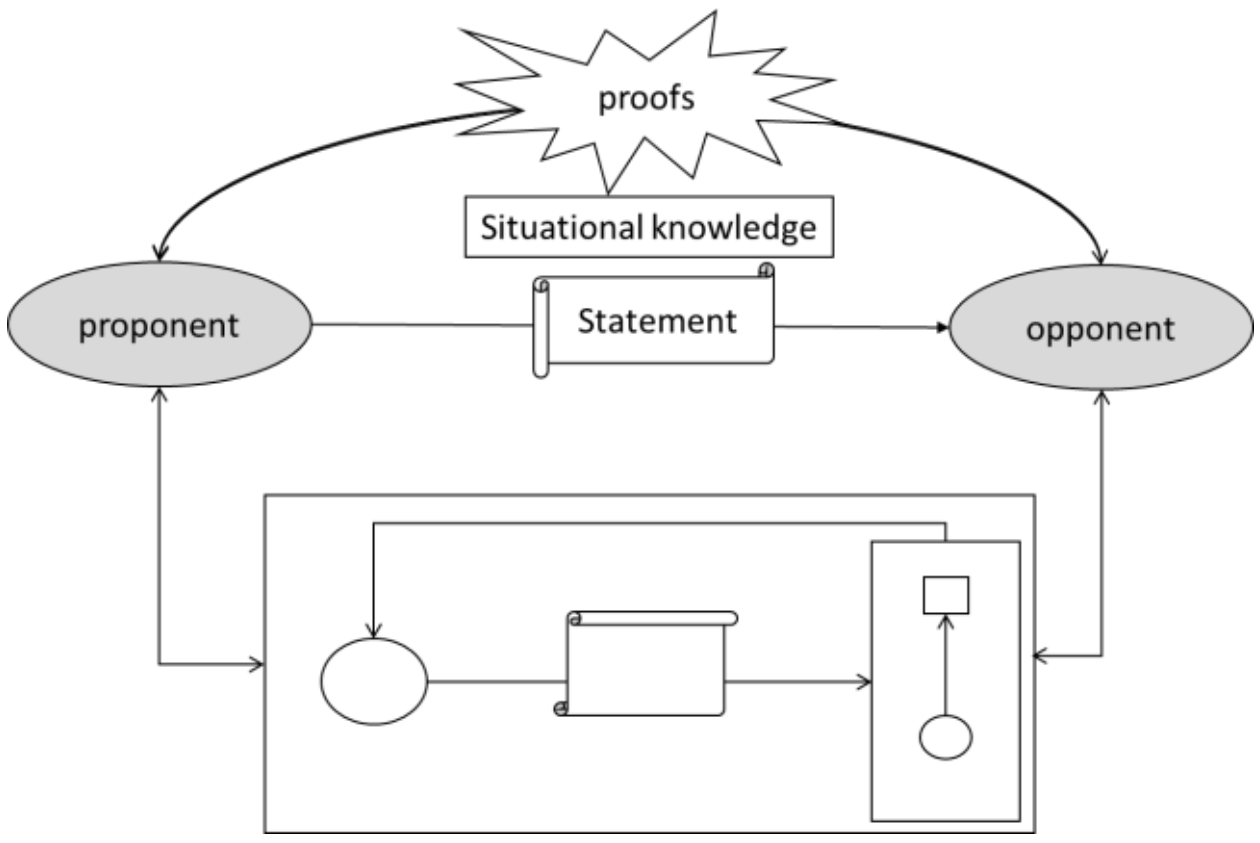

If an opponent who is providing feedback agrees with the statement, the proponent wins the game, but if the opponent has a valid objection, the opponent wins the game. These statements must be subject to judgment by the interlocutor, who cannot simply receive information. This means that the opponent must be able to provide feedback, and the proponent must be able to protest, reject reasons they judge to be false, and prove their position. Therefore, they must both be in symmetrical positions a priori: the validity of an objection must be approved by all the actants. The proponent and the opponent have the same access to information and action.

In section 5.2, inside a formulation situation, an additional demand is introduced for the receiver: to say if all information was useful or if some information was lacking, which represents a validation situation. In the validation situation schema (Fig. 5), the milieu is formed by the formulation situation, which includes the action situation. These are nested within the validation situation. The sentences that appeared during the formulation situation are the objects of the validation situation. The proof of the validation or rejection of the statement is shared by all actants. If the argument is validated at the end of the process, it becomes a theorem.

## 3.3. Two didactic processes

Now we revisit how the TS approaches adaptation and acculturation. As Brousseau has emphasized, the TS is not radical constructivism. However, he as well as numerous other researchers have dedicated significant time to an epistemological analysis of mathematics and mathematical situations (the TMS component) in order to comprehend all mathematics in terms of situations. In the TDSM, adaptation and acculturation are combined, and two complementary processes represent the necessary circulation between institutional knowledge (acculturation) and situational knowledge (adaptation). These processes are known as the institutionalization process and the devolution process.

### 3.3.1. Institutionalization process

As we saw above, several situations of different kinds are necessary to lead from a simple action (sometimes physical) on a milieu to an intellectual reflection as close as possible to mathematical reasoning. The teacher plays a crucial role, not only in choosing situations, but also in gradually making the most important findings, terms and concepts explicit and recording them in order to reach the knowledge to be taught. At the end of this institutionalization process, the teacher can formalize what is important, what needs to be retained, and how to reuse it. This may seem like a very traditional “lesson,” but because it comes at the end of a process based on the meaning of the mathematics involved, it's quite different: the outcome would not have the same meaning at all without the previously established situations.

### 3.3.2. Devolution process

When a teacher sets up a situation for the students to interact with a milieu, pursuing a given stake, regardless of the teacher's didactic intention, the teacher's actions are necessary not only at the beginning of the situation, but also as it unfolds. This is a delicate task, which is in some ways paradoxical (see 3.1, first paradox of teacher's decisions). From the students' point of view, it is not possible to distinguish between what is essentially adidactical (actions, formulations or validation in the milieu and stake setting) and what is of didactical origin (teaching intentions). "The final adidactical situation of reference, the one that characterizes the knowledge, can be studied in a theoretical way, but in the didactical situation, for the teacher as well as for the student, it is a sort of ideal towards which they are trying to converge" (Brousseau, 1997, p. 31). The adidactic situation is a theoretical model that enables us to characterize a given knowledge in terms of situation, not a description that may be followed step by step for teaching purposes.

As emphasized in the second section, the form of TS, as well as its didactic engineering methodology, is closely related to its project of understanding phenomena in the teaching and learning of mathematics. The next section and examples will provide further explanation.

# 4. Functions

The main function of the TS is to understand didactic phenomena, which must first be identified through observation. For this reason, most research using the TS framework employs a didactic engineering methodology. For about thirty years, the TS has also been used to study regular teaching and for teacher training. In the next section, we provide examples in geometry to illustrate the first aspect. Other examples may be found in Laborde, et al. (2005), for instance.

## 4.1. Didactic engineering to study didactic phenomena

The first objective of the TDSM is to identify, study and understand didactic phenomena, ultimately aiming to identify optimal teaching and learning conditions for each mathematical concept within the school curriculum. A didactic phenomenon involves the relationships between knowledge, a teacher who aims to transmit it and students who aim to acquire it. This includes studying the knowledge itself, its usefulness, and the problems it can solve. It also involves studying how this knowledge can be represented by a situation or a sequence of situations (see fundamental situation in section 2.3) in an accessible milieu for learners so as to construct adidactic situations within possible didactic situations. Finally, it involves studying how this representation confronts contingency through realizations in class. This methodology called "didactic engineering" (see 2.4) consists of four phases (Artigue, 2020):

1. Preliminary analyses of the mathematical subject at stake and of the learning conditions.

2. Creation of a didactic situation and a priori analysis of this situation with the help of theory.

3. Implementation and observation of this situation in classes.

4. A posteriori analysis of the gathered data, and confrontation with the a priori analysis to validate or modify elements of the milieu or the didactic organization.

The a priori analysis is a theoretical analysis of possible strategies for achieving the stake in the given milieu. It also predicts what procedures may be observed in class. The a posteriori analysis distinguishes between what is contingent and what is a characteristic of the milieu, given the students' prior knowledge. It assesses to what extent the milieu provides the necessary feedback for students to approach the new knowledge expected. The a posteriori analysis may also concern how the teacher deals with the two main processes: devolution and institutionalization. The didactic engineering methodology is described in detail in Section 5 using examples from geometry.

## 4.2. TDSM to study regular teaching and help teacher training

Later, the TDSM was also used to study, in regular teaching, teacher's practice and teacher's knowledge. In this case, the researchers' aim, therefore, is to understand what is happening in the situation being observed, rather than to create a didactic situation. First, observe a class. Then, based on your observations, identify the knowledge at stake, the milieu, and how this milieu might provide feedback to actions of the students. Building on Brousseau's work, Hersant and Perrin-Glorian (2005) studied teachers' practices by examining students' potential adidactic work according to the milieu's possible feedback, given students' supposed knowledge. They examined how the didactic contract evolves as the students' knowledge develops. They considered several dimensions of the contract, including the knowledge at stake (the mathematical domain and the didactic status of this knowledge for students) and the distribution of responsibility between the teacher and the students.

Margolinas et al. (2005) studied the situation from the teacher's perspective. They wondered how teachers could learn and acquire didactic knowledge by observing their students solve a problem. They used an extension of Brousseau's model of milieu structure to consider the teacher's different positions in relation to the milieu and the knowledge at stake in the situation.

# 5. Two examples of using TDSM in the context of geometry

Our aim in this section is to show an example of how TDSM can contribute to identify and to study didactic phenomena with a didactic engineering methodology. Geometry has many different functions in society and is interpreted in different ways within the school curriculum. Brousseau (2000b) distinguishes two main relationships to geometry: geometry as a model of space for practical purposes and geometry as a coherent mathematical theory based on an axiomatic foundation. To illustrate these, he imagines two different fundamental situations:

- A spatial situation: a carpenter working on the building of a cathedral must cut heavy pieces of wood on the ground that must be precisely adjusted for assembly ten meters above,
- A theoretical situation: when someone draws a triangle with a very obtuse angle and its three perpendicular bisectors, on paper using standard tools, due to the imprecision of the drawings, a very small triangle is usually formed. The challenge then becomes finding another triangle for which the triangle formed by the three bisectors is as large as possible.

The carpenter's situation represents useful knowledge for dealing with problems that arise and must be solved in the tangible space, even if it involves complex reasoning. In the case of the bisectors, searching for such a triangle leads to prove that it cannot exist: To ensure the model's logical coherence, the three bisectors must intersect at the same point. This situation is typical of theoretical geometry. Brousseau defines also three values for the variable size of space: micro-space (for instance drawings on a sheet of paper), meso-space: the space in which someone can move and control their environment by sight and the macro-space where no control by sight is possible.

Building on Brousseau's work, Berthelot and Salin (1998, 2000) distinguish between spatial knowledge, which is used to interact with space, and geometric knowledge. They also distinguish two approaches to the problems posed in the tangible space, particularly in the micro-space of a sheet of paper: the "practical problematics[3]" and the "modeling" or "spatio-geometric problematics". In the practical problematics, relations to space are effective, and controlled empirically and contingently by the senses. Validation is performed in the tangible space, often through perception. In the spatio-geometric problematics, the problem cannot be addressed directly in the tangible space: it must be translated into a geometric model and solved according to geometric criteria. Then, the result is translated back into the tangible space for validation. These two "problematics" differ from the geometric one in which the problem is presented directly within the context of geometry, where it is also solved, and the solution is validated within the framework of theoretical geometry, and according to established rules. Considering that the practical problematics hinders the spatio-geometric and the geometric problematics for students, Berthelot and Salin show that the size of the space is an important didactic variable, because controlling by sight or with material instruments is not equally efficient depending on the size of the space. For example, half of the observed students aged 10-11 were unable to predict the position of the feet of a rectangular bench after it had been moved, even though most of them could draw a rectangle on a sheet of paper, as visual control can be used continuously in this case (Berthelot & Salin, 1998, p. 76-77). Berthelot and Salin propose the spatio-geometric problematics as a way to enter the geometric one through problems posed in space, while blocking natural spatial procedures, especially by working in the meso-space. They thus retain three important classes of variables for analyzing geometric situations in education: the type of problematics, the type of relationship to space, effective or internalized, the adidactic character or not of situations.

The difficulties involved in learning geometry, including how geometric drawings should be viewed in theoretical geometry, have been studied for a long time. From a cognitive or epistemological point of view, Van Hiele identified five levels to describe how students reason about shapes and geometric ideas. Fischbein (1993) introduced the construct of figural concept, and Houdement and Kuzniak (2006) spoke of different paradigms. Duval (2014) identified the articulation of semiotic registers in geometry (graphic / verbal register) as a major difficulty. Research in which one of the authors of this chapter participated (Mathé et al., 2020; Mathé & Perrin-Glorian, 2023; Perrin-Glorian et al., 2024) focused on the following question: how can geometry teaching be managed in primary and first grades of secondary schools to help students' transition to theoretical geometry? This research produced situations that allowed us to address certain aspects of geometry learning that are insufficiently covered in usual teaching. Below, we present the didactic engineering approach used in this research, providing a more specific example for classroom implementation. We will supplement this with an example from secondary education.

[3] Here, "problematics" refers to a set of problems and the way they are expected to be approached.

## 5.1. Transition to theoretical geometry

### 5.1.1. Preliminary analyses

Preliminary analyses focused mainly on two areas: firstly, an epistemological reflection on geometry itself and how to organize a coherent distribution of knowledge for students throughout the curriculum according to their age (Karp & Perrin-Glorian, to appear); and secondly, a review of previous studies on geometry learning. Next, we elaborated on the main points that emerged from these analyses.

There is a gap between reproducing or studying geometric drawings with geometric tools and theoretical geometry, concerning objects and how to validate affirmations. Nevertheless, many properties are expressed in the same way because the same words are used for different objects. Language and drawings on a two-dimensional plane play an essential role in geometry. Drawings can represent space, material objects or geometric shapes. They can be created using geometric tools or dynamic geometric environments. In theoretical geometry, both verbal and graphic registers are required: in proofs, drawings can serve various functions, such as exploring multiple configurations for conjectures or adding new lines to reveal relationships or isolating sub-figures. A major difficulty is that the visualization of figures in geometry differs from the usual way of looking at and talking about drawings. Duval (2014) shows that, in a geometric activity, the relationship between visualization and discourse requires the dimensional deconstruction of shapes: a drawing must generally be seen as composed of different other shapes, a 2D-shape must be perceived as a network of lines (1D) and points (0D). Visualization of figures involves cognitive processes, but it is also closely related to geometric knowledge and requires coordination between graphic and verbal registers. Indeed, “the mathematical way of seeing” drawings utilizes complementary features expressed in words (Mathé & Mithalal, 2019; Perrin-Glorian et al., 2024).

These preliminary analyses led to the hypothesis that certain aspects of geometry learning should be addressed more explicitly in teaching. For example, the necessity to see a point as the intersection of two lines (straight lines or circles) or to see a straight line under a segment. We took this issue as a main stake to create a milieu.

### 5.1.2. Construction of a milieu

The approach to geometry adopted for students was based on the "micro-space" of the sheet of paper and the "milieu of plane figures" (Fregona, 1995). A fundamental situation was defined by identifying didactic variables in the problem of reproducing plane figures, mainly the available tools and the reproduction conditions. These didactic variables were chosen so as to block the practical problematics and orient the students towards the spatio-geometric problematics. In addition, this research hypothesizes that to approach geometric concepts, reproduction must be produced without measuring. On the one hand, it dismisses possible difficulties with numbers; on the other hand, when numbers are involved, they capture the students' full attention.

To represent the basic concepts of Euclidean plane geometry without measurement, “theoretical instruments” were defined (Perrin-Glorian, 2024). These instruments are theoretical in that they do not have the limitations of standard material tools (for example, as a theoretical instrument, the ruler is as long as desired). These theoretical instruments are essential for the mathematical part of the a priori analysis. They include a ruler to draw straight lines, a compass to draw circles, an instrument to copy lengths, and another instrument to copy angles. These theoretical tools can be physically produced in several ways for use in the classroom.

An additional didactic variable is obtained by providing students with a pre-drawn part of the figure to be reproduced which may or may not be the same size as the model. Students must reproduce the model from this "starting figure." In addition to the available tools, the didactic variables concern the geometric properties of the figure to be reproduced (model) and the choice of the starting figure. Both depend on the knowledge at stake.

To obtain an adidactic situation of action, in which the milieu, rather than the teacher, controls the action, a copy of the model figure on transparent paper is only made available to the students when they believe the reproduction is complete; the copy is kept away from the students. The game consists of reproducing the model from the starting figure with the tools. You win the game if the figure on the transparent paper exactly fits with your reproduction. To obtain an adidactic situation of formulation, a student (the sender) must provide information to another student (the receiver) so that the receiver can reproduce the figure that they cannot see. Then the sender and the receiver compare their figures. During this comparison, the interaction with the milieu is effective but, to write a message, the sender must refer to shared knowledge that the receiver can understand. Likewise, the receiver must refer to this knowledge to understand the message. They must use a shared language.

### 5.1.3. Implementation in classrooms

Examples of implementations of such situations in primary classes may be found in Mangiante-Orsola (2023); Mathé (2025); Mathé et al. (2020); Mathé & Perrin-Glorian (2023); Perrin-Glorian et al. (2024). We present here a sequence of four situations with students 9-12 years old about conceptualizing lines and points (The main objective of this sequence is that students learn to identify and understand the role of lines and points in defining a geometric figure, as well as the relationships between them. The sequence is organized into four situations, each with a different starting figure and different tools, in addition to the ungraduated ruler that is always available. The milieu is constituted by the model figure, the starting figure and the tools available. The stake is always to reproduce the same model figure: a quadrilateral with no special properties, from which triangles are constructed using its diagonals (Fig. 6). To facilitate the written descriptions and analyses, a copy with names for points is included here (Fig. 7).

*Figure 6. Model figure for students*

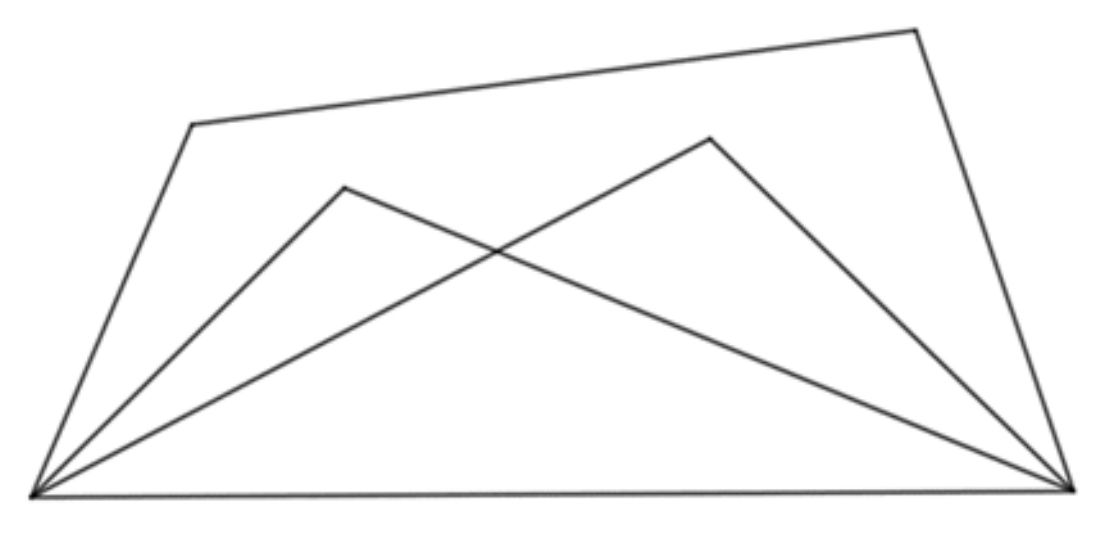

*Figure 7. Model figure, not available for students*

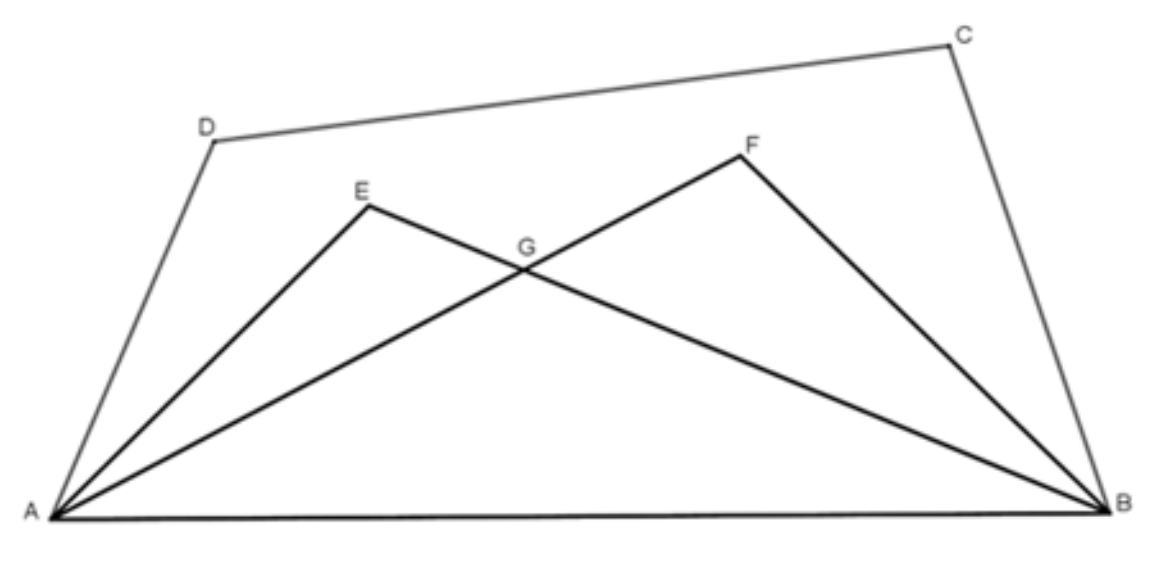

The milieu changes according to the rules of the game, depending on the starting figure and the available tools. This results in different learning objectives. At the end of each situation, there is a collective moment for sharing solutions with local institutionalization of some results. These moments introduce shared knowledge that the students can then draw on in their interactions with the milieu.

The *first situation* aims to understand the model figure. The starting figure is the quadrilateral ABCD only. The instruments provided are an ungraduated ruler and the AEB and AFB triangle templates (Fig. 8). The stake is to identify additional triangles beyond the three that appear juxtaposed at first glance in the model figure.

*Figure 8. The starting figure and the templates.*

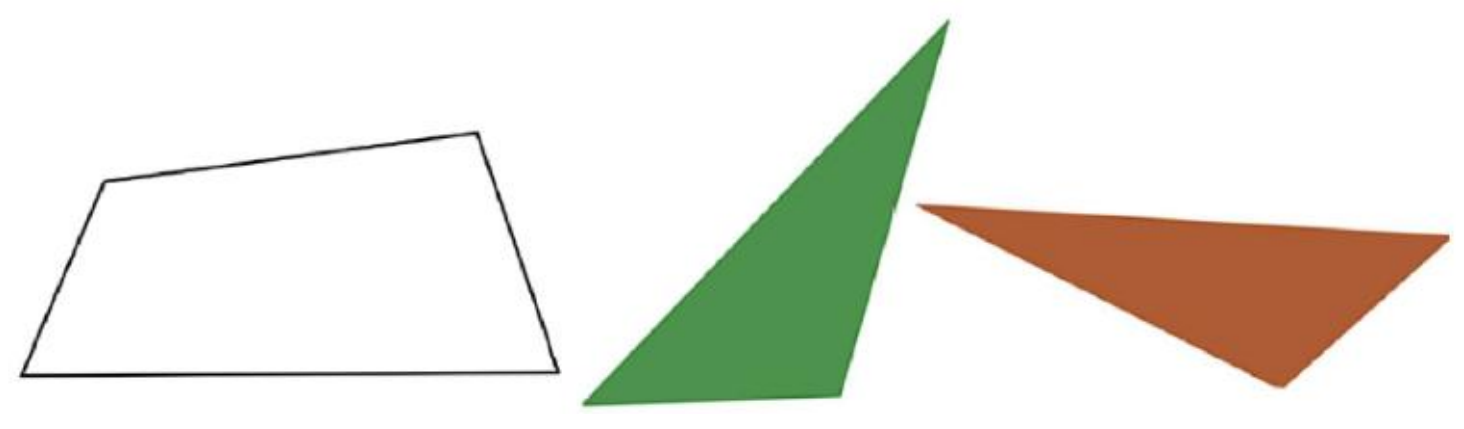

*Figure 9. How to place the templates.*

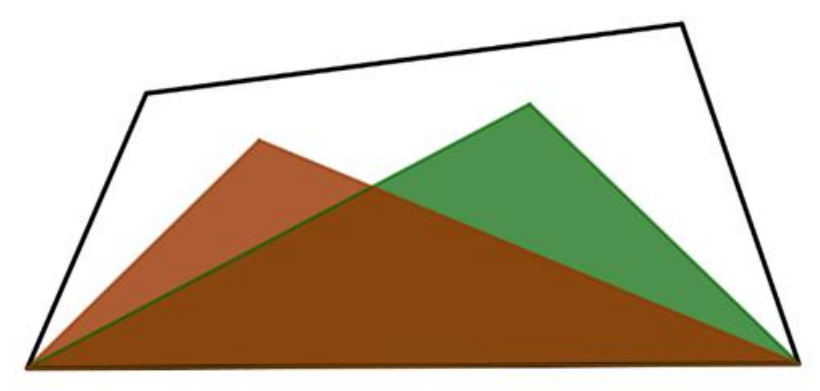

To succeed, students must place each template on the model figure to identify the corresponding triangle. Then, they must place the templates on the starting figure to reproduce the triangles of the model (Fig. 9). The difficulty of this strategy lies in the necessity of overlapping the templates, which conflicts with the usual didactic contract of using templates in puzzles. Thus, this situation also aims to establish a new didactic contract with two ruptures. Students are accustomed to using a graduated ruler to reproduce figures by measuring lengths but now, they must use an ungraduated ruler and templates. Moreover, the templates must overlap. It's a big step for students. We anticipate that the teacher may need to encourage them and explicitly allow overlap. However, the teacher must continue to reject the graduated scale because it is necessary to maintain the milieu. The teacher plays an important role in helping students accept these changes, this role is part of the devolution process. The situation concludes with a local institutionalization: In the model figure, we can see a quadrilateral with five triangles inside.

The *second situation* aims to identify the alignments. In the starting figure of the new milieu, the triangle T1 is already drawn, and the tools are the ungraduated ruler and a template of T2 with two corners cut out (Fig. 10).

*Figure 10. Starting figure and template.*

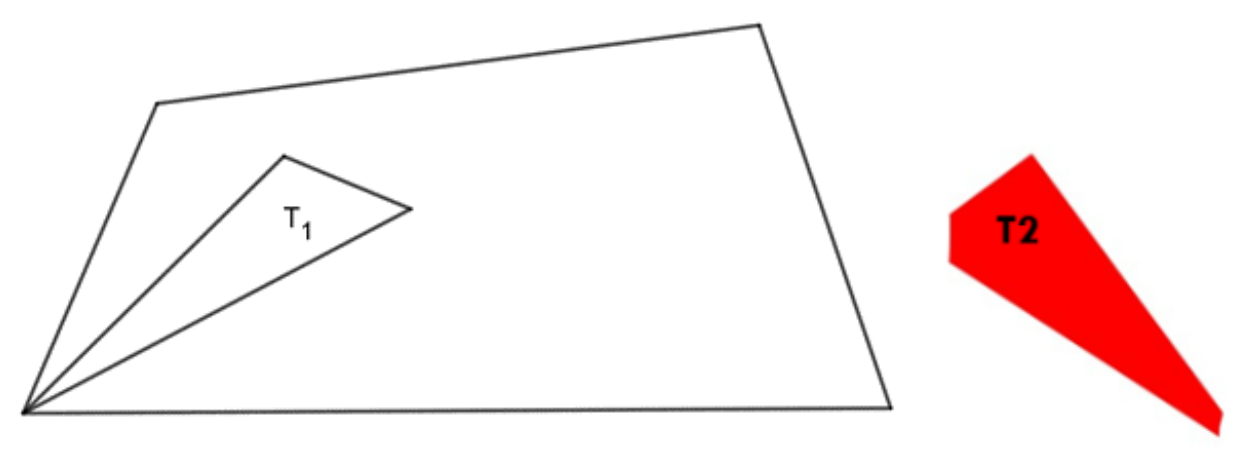

Now, the objective is to identify the alignments of A, G, F, on one part and E, G, B on the other. In the first situation, these alignments were supported by the milieu as sides of the large triangles. Now, they must be identified to extend the sides of T1 and set the template in the correct position before tracing T2 onto the starting figure. The corners of the template are cut to prevent adjusting it by coinciding one vertex of T2 with the vertex of T1 at G and another with the vertex of the quadrilateral at B. This makes it necessary to use a ruler for extending the lines. The concept of angle as defined by its sides, even without a vertex, is implicitly at stake though it is not expected to institutionalize it. Some students can extend the triangle's sides up to the quadrilateral's vertices, but it is not the focus here. The situation concludes with a local institutionalization: we need to extend the sides of the triangle to place the template and draw the second triangle.

The objective of the *third situation* is to understand that a point must not be placed on a line by sight alone. It must be obtained by transferring a length or by intersecting another line. In the new starting figure, the top of the quadrilateral is erased (Fig. 11) and this figure is drawn at a larger scale than the model to prevent students from transferring lengths from the model to the starting figure. The only tool available is an ungraduated ruler.

To draw the missing side of the quadrilateral, the students must construct its two extremities (Fig. 12). Since the model and the starting figure have different scales, it is impossible to transfer lengths. The only way to succeed is to identify new alignments: A, G, F, and C (respectively B, G, E, and D). Then, the students can find C (or D) by intersecting an extension of the visible part of a side of the quadrilateral and the extension of [AF] (or [BE]).

*Figure 11. Starting figure.*

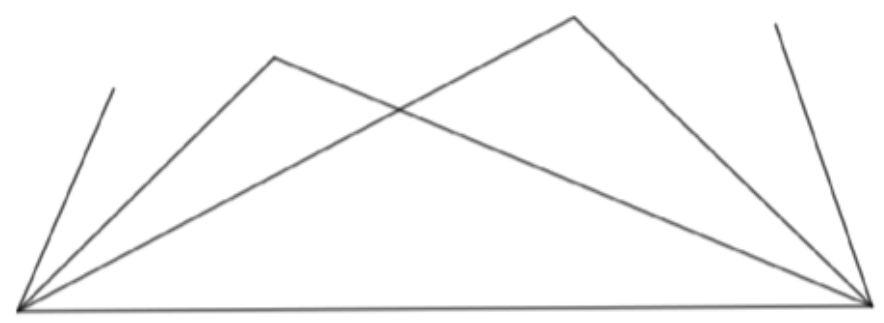

*Figure 12. Lines to draw.*

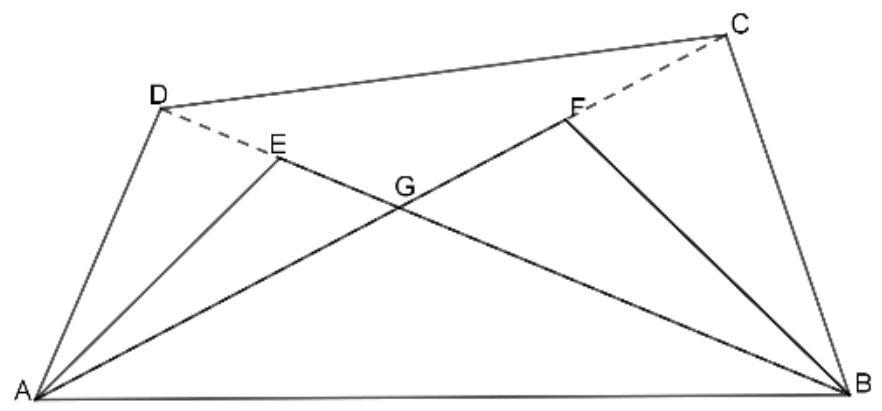

This is the crucial part of the sequence for students in grades four through six. We anticipate that most students will extend the incomplete sides of the quadrilateral and place C and D on them, judging their position by sight. However, checking with the model's copy on transparent paper shows that it does not fit. We anticipate that some students will try to transfer the length of [BC] with the edge of a sheet of paper and that most of them will have difficulty finding another line that intersects the first one if they have not extended the line [AF] (or line [BE]) up to point C (or D) in the second situation. Thus, the teacher will have to remind students that only the ungraduated ruler may be used and to encourage them, which is part of the devolution process.

This situation yields three conclusions regarding local institutionalization: two sides of the small triangles lie on the quadrilateral's diagonals; the sides of the triangles can be extended as long as desired to form a straight line; and to determine the missing vertex, a specific length is required. If this length cannot be obtained, one should look for another line on the model that intersects the first line at the vertex. A point can be found by intersecting two lines.

In the *last situation* of the sequence, no new knowledge is at stake, this situation is part of the institutionalization process. This situation provides an opportunity for the teacher to ensure that the knowledge has effectively disseminated throughout the class and that all students understand the role of the quadrilateral's diagonals in constructing the figure. The new starting figure is only the quadrilateral ABCD, the tools are the ungraduated ruler and templates of triangles T1 and T2 with cut corners.

In the previous situation, the quadrilateral's diagonals were identified, so we anticipate that most students will draw them. Placing the templates into an angle of these diagonals is a reinvestment of the second situation. Then completing the triangles will be easy.

The final part of the institutionalization process aims to decontextualize the knowledge at stake, in an interaction between the teacher and the students: A line can be extended indefinitely, as long as desired; to find the second endpoint of a segment on a line, we need a length or another line to intersect the first one. A point is formed by intersecting two lines.

## 5.1.4. A posteriori analysis

The sequence presented in 5.1.3. was observed in several classes and video recorded in some of them. In the a posteriori analyses of these situations, we compared students' actual procedures with those predicted in the a priori analysis. We also analyzed the teacher's interventions to identify which led to a change of the initial milieu. The sequence of four situations presented above is the outcome of observations of several implementations with slight variations coming from a posteriori analysis of previous versions. As expected, teachers must intervene in situation 1 to allow the overlapping of templates. The situation 2 presented above results from a first implementation where the template was intact. Many students adjusted the corners so that extending the triangle's sides was not needed. It is why we decided to cut the corners. In situation 2 as described above, not all students immediately extended the lines. Some students extended the segment [AG] and placed the template's edge on the extension. However, they placed the cut corner in G, so the second edge was under the diagonal and not in the correct position. They could realize this by placing the model traced on transparent paper over their figure.  The situation 3 presented above also resulted from an initial implementation in which the starting figure was different. Some students need help with the current version. The teacher must encourage them to identify other helpful lines on the model.

In our example, the objective was to create teaching resources. Thus, we sought to determine whether the situation enabled students to learn the targeted knowledge. However, didactic engineering can be used to answer other research questions. For instance, in the above-mentioned research by Berthelot and Salin, the question was: Are students able to use their scholarly geometric knowledge to address spatial problems? The bench situation revealed that half of the students could not.

## 5.2. Congruence of triangles. Formulation and validation situations

This sequence was designed to help eighth-grade students (ages 13–14) understand that a triangle can be defined by three specific measurements (Perrin-Glorian & Didier, 2022). The objective is for students to identify the three cases of triangle congruence so that they can use them in proofs later. The adidactic situation is a combination of formulation and validation situations. Students are given information about a triangle, including the measurements of its three sides and three angles. First, they must construct the triangle and then write a message with as little information as possible so that another student can construct the same triangle. Finally, they compare the two triangles. Moreover, the receivers must say if they used all the information to construct their triangle. For each student, there is a part of action (to draw a triangle from information on measurement). The formulation situation consists of the sender selecting the information to be conveyed to the receiver and then verifying with the receiver that the triangles constructed by the sender and the receiver can be superimposed. The additional questions give as little information as possible for the sender and say if all information was useful or if some information was lacking, for the receiver, opens the validation situation that leads to a debate in the classroom afterwards.

For a teacher in a classroom with about thirty students, organizing this formulation situation is very demanding. The teacher must select the triangles and manage the exchange of messages. Students exchanging messages must have different triangles to describe. Furthermore, it is preferable that they make different choices when describing their triangles. Such a situation is organized when the knowledge at stake is significant as well as its explicit formulation. Selecting properties that are necessary and sufficient to characterize a triangle is an important challenge for understanding cases of triangles congruence as tools for proof.

### 5.2.1. A priori analysis and choice of didactic variables

Students are supposed to be able to construct triangles with a graduate ruler, a compass and a protractor (prior knowledge). The stake is to write the requested message and to agree with the receiver that this message was as efficient as possible. We hope the students' messages cover the three cases of triangles congruence with examples of non-sufficient information or non-necessary information. We also hope to find a case with two length measurements and one angle measurement where the angle is incorrectly chosen (i.e. not between the measured sides), resulting in two triangles that satisfy the conditions.

The number and selection of triangles are didactic variables in this situation. In the observed classroom, the teacher selected four triangles and organized pairs of students, which ensured sufficient diversity. Furthermore, all triangles had three acute angles (Fig. 13), which increased the likelihood of encountering a case with two solutions. Indeed, if a triangle has an obtuse angle, then any two length measurements and one angle measurement determine that triangle. Students were given a drawing with a triangle with measurements indicated but the drawing did not fit with them (Fig. 14).

*Figure 13. Measurements of the sides and angles of the four triangles.*

| $AB$ | $BC$ | $AC$ | [illegible] | [illegible] | [illegible] |
|---|---|---|---|---|---|
| 6,8 | 8 | 7,1 | 70° | 57° | 53° |
| 5,6 | 7,4 | 6,1 | 78° | 54° | 48° |
| 6,8 | 10,7 | 10,1 | 76° | 66° | 38° |
| 5 | 7,1 | 5,8 | 82° | 54° | 44° |

*Figure 14. Example of instructions for the sender.*

Write a text with as little information as possible so that the reader can construct the triangle.

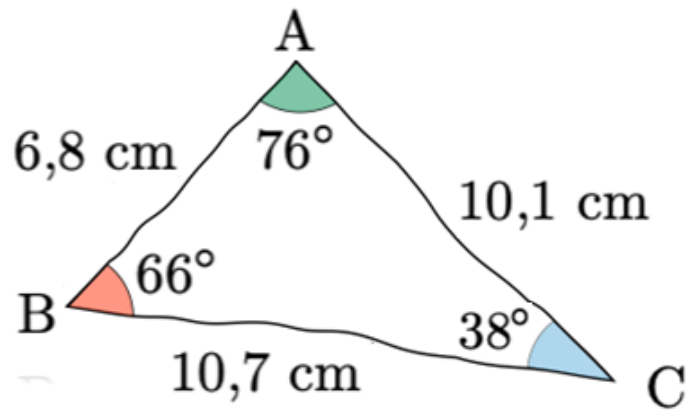

The milieu depends on the actant in question, and it evolves. For the sender, it consists of a triangle with information on measurements of all lengths and angles (Fig. 14). For the receiver, it consists of information on a triangle to draw but it's not clear if these pieces of information are sufficient or necessary. For the comparison, the milieu consists of the two triangles constructed, the one by the sender, the other by the receiver. For the discussion, inside the validation situation, it consists of all previous milieus and the receiver's statement on usefulness of information given.

To organize the exchanges so that each pair would receive a different message and have a different set of information, the teacher had the students write the messages during one class session and exchange them during the next one. He prepared a grid for analyzing students' messages between the two sessions.

### 5.2.2 Validation situation and institutionalized statements

All cases anticipated in the a priori analysis were encountered. The teacher first gave the floor to the pairs who had provided too much or too little data, then to those who had written a message corresponding to a case of congruence. A case with two solutions was also encountered. The conclusion of the formulation situation was that three pieces of information are sufficient to construct a triangle, provided they are chosen carefully.

Then, the teacher encouraged the students to formulate more general statements that can be written in their lesson notebooks. They engaged the students in a validation situation, involving a debate in which counterexamples were brought up, particularly regarding the choice of angles.

Finally, the class arrived at the following statement: Two triangles can be superimposed if either the lengths of two sides and the angle formed by these sides are equal, or two angles and the length of the segment joining the vertices of these angles are equal, or the lengths of all three sides are equal.

### 5.2.3 Institutionalization process: reinvestment in another proof

The students encountered cases of congruence when constructing triangles. However, other situations are needed for them to become a tool for proof. Here is an example (Fig. 15) that was used shortly afterwards in the same class. The students knew that the sum of the angles of a triangle is 180°.

*Figure 15. Use the result to prove.*

In the following figure, calculate the perimeter of quadrilateral ABCD.
Justify your answer.

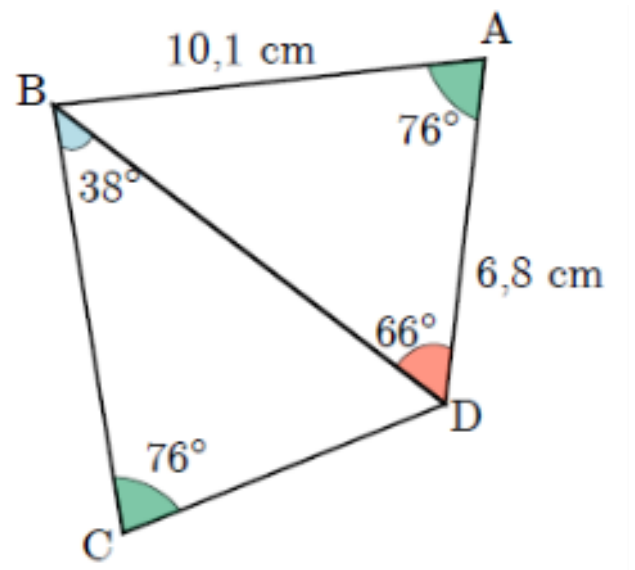

# 6. Conclusion

Guy Brousseau laid the first foundations of the TS in the 1960s and 1970s. He developed and stabilized them in the 1980s and 1990s in the context of interactions with students, teachers, educators, and researchers who visited the COREM or invited him and his wife Nadine Brousseau to visit other countries and discover other educational research cultures. Since then, many researchers have contributed to the development of various aspects of TS. However, the main concepts have not changed significantly over the years because they provided a solid foundation and a coherent way to observe and analyze didactic phenomena.

The TS, with its two branches, the TMS and the TDSM, provides now a robust network of concepts, methods, and results that help us better understand and predict what might happen when teaching mathematics at school. We hope we succeeded in showing why it constitutes a theory.

The TMS has produced many results that have been widely used in subsequent research in France and beyond. Epistemological-didactic analyses have been conducted in various fields to identify situations representing the core of mathematical concepts, such as whole numbers, rational and decimal numbers, functions, and probabilities. Some TMS studies have also identified transparent knowledge, such as enumeration, spatial knowledge, and reasoning (Margolinas & Laparra, 2011). The TDS based on the didactic engineering methodology has produced a wide range of situations implemented in controlled conditions. Even though these situations are not intended to be directly replicated in regular teaching, nevertheless, these situations can serve as a resource for teacher training and introduce certain concepts to help teachers reflect on their own practice (Mangiante-Orsola et al., 2018).

The TS has been used in numerous research studies —sometimes in combination with other theoretical frameworks or in disciplines outside of mathematics— for a very long time in France and around the world for decades. It continues to be developed and attempts of networking the TS with other theories in international groups provided fruitful research (Bikner-Ahsbahs & Prediger, 2014, chapters 4, 10, 12).

We really believe the best way to understand the TS is to view it through its functions, which are to study and understand didactic phenomena, in an experimental setting or in regular teaching. The TS provides the means to create situations in which these phenomena can be observed in rich environments. This is why we devoted a large part of this chapter to developing an example.

In an attempt to outline the form of the theory, we drew mainly upon the most well-known concepts. We hope the example will help readers recognize the underlying connections between these concepts. Although studying the TS, understanding it, and taking advantage of its potential can be difficult for beginners, we hope it will provide them with a new approach to learning and teaching mathematics. A course created at the university of Grenoble that was recently translated into English may be helpful (Bessot, 2003), as well as the ICMI Amor unit concerning Brousseau (https://www.mathunion.org/icmi/awards/amor/awardees-units).

# 7. References

Artigue, M. (2020). Didactic engineering in mathematics education. In S. Lerman (Ed.), *Encyclopedia of Mathematics Education* (pp. 202–206). Springer. https://doi.org/10.1007/978-3-030-15789-0_44

Berthelot, R., & Salin, M.-H. (1998). The role of pupils' spatial knowledge in the elementary teaching of geometry. In C. Mammana & C. Villani (Eds.), *Perspective on the teaching of geometry for the 21st century. An ICMI study* (pp. 71–78). Kluwer Academic Publishers. https://doi.org/10.1007/978-94-011-5226-6_3

Berthelot, R., & Salin, M.-H. (2000). L'enseignement de la géométrie au début du collège. Comment concevoir le passage de la géométrie du constat à la géométrie déductive ? [Teaching geometry at the start of secondary school. How should we approach the transition from observational geometry to deductive geometry?] *Petit x*, *56*, 5–34.

Bessot, A. (2003). Une introduction à la théorie des situations didactiques [Introduction to the theory of didactical situations]. *Les cahiers du laboratoire Leibniz*, *91*. https://telearn.archives-ouvertes.fr/hal-00078794/document available in English at https://hal.science/hal-04500947v1

Bessot, A., & Margolinas, C. (2025). Situações fundamentais no centro da teoria das situações matemáticas [The fundamental situations at the core of the theory of mathematical situations.]. In S. Ag Almouloud, P. Job, L. M. Santos Farias, & A. Henriques (Eds.), *Teoria das situações: Fundamentos teóricos e metodológicos para a investigação e a formação profissional* [Theory of situations: Theoretical and methodological foundations for research and professional training] (pp. 89–120). Editora CRV. Version française https://hal.science/hal-05575541

Bikner-Ahsbahs, A., & Prediger, S. (Eds) (2014). *Networking of Theories as a Research Practice in Mathematics Education*. Springer. https://doi.org/10.1007/978-3-319-05389-9

Bloch, I. (2002). Différents niveaux de modèles du milieu dans la théorie des situations [Different levels of milieu models in the theory of situations]. In J.-L. Dorier, M. Artaud, M. Artigue, R. Berthelot, & R. Floris (Eds.), *Actes de la 11ème école d'été de didactique des mathématiques* [Proceedings of the 11th summer school on the didactics of mathematics] (pp. 125–140). La Pensée Sauvage.

Brousseau, G. (1980). L'échec et le contrat [Failure and the contract]. *Recherches*, *41*, 177–182.

Brousseau, G. (1997). *Theory of didactical situations in mathematics* (N. Balacheff, M. Cooper, R. Sutherland, & V. Warfield, Eds.). Kluwer Academic Publishers.

Brousseau, G. (2000a). Educacion y didactica de las matematicas [Education and didactics of mathematics] (D. Block & Martinez Falcon, Trans.). *Educacion matematica*, *12*(1), 5–38. Version française https://hal.archives-ouvertes.fr/hal-00466260/fr/

Brousseau, G. (2000b). *Les propriétés didactiques de la géométrie* élémentaire [The didactic properties of elementary geometry]. Rethymnon, Crete. https://hal.science/hal-00515110

Brousseau, G. (2004). Felix Klein Medallist: Research in mathematics education. In M. Niss (Ed.), *Proceedings of the Tenth International Congress on Mathematical Education* (pp. 244–254). IMFUFA, Roskilde University.

Brousseau, G. (2006). Mathematics, didactical engineering and observation. In J. Novotná, H. Moraová, M. Krátká, & N. Stehlíková (Eds.), *Proceedings of the 30th Conference of the International Group for the Psychology of Mathematics Education* (Vol. 1, pp. 1.3–1.18). online IGPME. https://www.igpme.org/publications/current-proceedings/

Brousseau, G., Brousseau, N., & Warfield, G. (2014). *Teaching Fractions through Situations: A Fundamental Experiment*. Springer. https://doi.org/10.1007/978-94-007-2715-1

Derouet, C., Durand-Guerrier, V., & Mathé, A.-C. (Eds) (2025). *Pré-actes du Colloque international en hommage à l'œuvre de Guy Brousseau* [Preliminary proceedings of the international symposium in honour of the work of Guy Brousseau]. IREM d'Aquitaine. https://hal.science/hal-05184913v1

Duval, R. (2014). The first crucial point in geometry learning: Visualization. *Mediterranean Journal for Research in Mathematics Education*, *13*(1-2), 1–28.

Fischbein, E. (1993). The theory of figural concepts. *Educational Studies in Mathematics*, *24*(2), 139–162. https://doi.org/10.1007/bf01273689

Fregona, D. (1995). *Les figures planes comme « milieu » dans l'enseignement de la géométrie : Interactions, contrats et transpositions didactiques* [Plane figures as a 'milieu' in the teaching of geometry: Interactions, contracts and didactic transpositions]. Université de Bordeaux I.

Gascòn, J., & Nicolás, P. (2017). Can didactics say how to teach? The beginning of a dialogue between the anthropological theory of the didactic and other approaches. *For the Learning of Mathematics*, *37*(3), 9–13.

Hersant, M., & Perrin-Glorian, M.-J. (2005). Characterization of an ordinary teaching practice with the help of the theory of didactic situations. *Educational Studies in Mathematics Education*, *59*, 113–151. https://doi.org/10.1007/s10649-005-2183-z

Houdement, C., & Kuzniak, A. (2006). Paradigmes géométriques et enseignement de la géométrie [Geometric paradigms and the teaching of geometry]. *Annales de Didactique et de Sciences Cognitives*, *11*, 175–195.

Karp, A., & Perrin-Glorian, M.-J. (to appear). Curricular ideas and changes: different times, different countries, and different grades. In A. Gutiérrez & T. Lowrie (Eds), *26th ICMI study. Advances in geometry education, section 5, Ch. 22*.

Laborde, C., Perrin-Glorian, M.-J., & Sierpinska, A. (Eds.). (2005). *Beyond the Apparent Banality of the Mathematics Classroom*. Springer. https://doi.org/10.1007/0-387-30451-7

Mangiante-Orsola, C. (2023). Study on teacher appropriation of a geometry education resource. In C. Guille-Biel Winder & T. Assude (Eds.), *Articulations between sensitive space, graphic space, geometric space. Resources, practices and training* (pp. 199–219). ISTE Science Publishing. https://doi.org/10.1002/9781394229604.ch11

Mangiante-Orsola, C., Perrin-Glorian, M.-J., & Strømskag, H. (2018). Theory of didactical situations as a tool to understand and develop mathematics teaching practices. *Annales de Didactique et de Sciences Cognitives. Special issue English-French,* 145–174. https://doi.org/10.4000/adsc.334

Margolinas, C. (2002). Situations, milieux, connaissances : Analyse de l'activité du professeur [Situations, contexts, knowledge: An analysis of the teacher's work]. In J.-L. Dorier, M. Artaud, M. Artigue, R. Berthelot, & R. Floris (Eds.), *Actes de la 11ème Ecole d'Eté de*

*Didactique des Mathématiques* [Proceedings of the 11th summer school on the didactics of mathematics] (pp. 141–156). La Pensée Sauvage. https://hal.archives-ouvertes.fr/halshs-00421848

Margolinas, C. (2022). Theory of didactic situation in mathematics: An epistemological revolution. In Y. Chevallard, B. Barquero, F. Ignasi, J. Gascòn, P. Nicolás, & N. Ruiz-Munzón (Eds). *Advances in the Anthropological Theory of the Didactic* (pp. 37–49). Springer Nature. https://hal.science/hal-03606376 https://doi.org/10.1007/978-3-030-76791-4_4

Margolinas, C., & Bessot, A. (2024). Modules 1 to 11. In C. Margolinas & A. Bessot (Eds.), *ICMI AMOR Project – Guy Brousseau Unit – The Theory of Didactical Situations in Mathematics*. https://www.mathunion.org/icmi/awards/amor/guy-brousseau-unit

Margolinas, C., & Bessot, A. (2025). A estruturação do milieu: Metodologia de análise de situações didáticas [The structuring of the milieu: A methodology for analysing didactic situations]. In S. Ag Almouloud, P. Job, L. M. Santos Farias, & A. Henriques (Eds.), *Teoria das situações: Fundamentos teóricos e metodológicos para a investigação e a formação profissional* [Theory of situations: Theoretical and methodological foundations for research and professional training] (pp. 121–154). Editora CRV. Version française https://hal.science/hal-04749448

Margolinas, C., Coulange, L., & Bessot, A. (2005). What can the teacher learn in the classroom? *Educational Studies in Mathematics*, *59*(1), 205–234. https://doi.org/10.1007/s10649-005-3135-3

Margolinas, C., & Laparra, M. (2011). Des savoirs transparents dans le travail des professeurs à l'école primaire [Transparent knowledge in the work of primary school teachers]. In J.-Y. Rochex & J. Crinon (Eds.), *La construction des inégalités scolaires* [The construction of school inequalities] (pp. 19–32). Presses Universitaires de Rennes.

Mathé, A.-C. (2025). *Enseigner et apprendre la géométrie à l'École. Itinéraire et perspectives d'une recherche en didactique des mathématiques* [Teaching and learning geometry at school. A journey and prospects of research in mathematics didactics] Université Clermont-Auvergne. https://hal.science/tel-05052009v1

Mathé, A.-C., Barrier, T., & Perrin-Glorian, M.-J. (2020). *Enseigner la géométrie à l'école élémentaire—Enjeux, ruptures et continuités* [Teaching geometry in primary school—Challenges, changes and continuities]. Academia - L'Harmattan.

Mathé, A.-C., & Mithalal, J. (2019). L'usage des dessins et le rôle du langage en géométrie : Quelques enjeux pour l'enseignement [The use of drawings and the role of language in geometry: Some challenges for teaching]. In S. Coppé, É. Roditi, V. Celi, F. Chellougui, F. Tempier, C. Allard, C. Torriveau, M. Haspekian, P. Masselot, S. Rousse, H. Sabra, & M. Tiwan-Zaka (Eds.), *Nouvelles perspectives en didactique : Géométrie, évaluation des apprentissages mathématiques* [New Perspectives in Teaching and Learning: Geometry, Assessment of Mathematical Learning] (pp. 47–85). La Pensée Sauvage.

Mathé, A.-C., & Perrin-Glorian, M.-J. (2023). The geometry of tracing, a possible link between geometric drawing and Euclid's geometry? In C. Guille-Biel Winder & T. Assude (Eds.), *Articulations between sensitive space, graphic space, geometric space. Resources, practices and training* (pp. 3–33). ISTE Science Publishing. https://doi.org/10.1002/9781394229604.ch1

Niss, M., & Jankvist, U. T. (2022). On the Mathematical Competencies Framework and Its Potentials for Connecting with Other Theoretical Perspectives. In U. T. Jankvist & E. Geraniou (Eds). *Mathematical Competencies in the Digital Era* (pp. 15–38). Springer. https://doi.org/10.1007/978-3-031-10141-0_2

Perrin-Glorian, M.-J. (1994). Théorie des situations didactiques : naissance, développement, perspectives [Theory of didactic situations: origins, development, prospects] in M. Artigue, R. Gras, C. Laborde, P. Tavignot (Eds). *Vingt ans de didactique des* *mathématiques en France* [Twenty years of mathematics didactics in France], pp.97–147. La Pensée Sauvage.

Perrin-Glorian, M.-J. (2024). Teaching plane geometry in coherence from 6 to 15 years. *Rematec*, *19*(48). https://doi.org/https://doi.org/10.37084/REMATEC.1980-3141.2024.n48.e2024001.id588

Perrin-Glorian, M.-J., & Didier, G. (2022). Une organisation de la géométrie du collège. Les cas d’isométrie et de similitude des triangles comme outils de démonstration [An approach to teaching geometry in secondary school. Isometry and similarity of triangles as tools for proof]. *Petit x*, *117*, 21–51.

Perrin-Glorian, M.-J., Mathé, A.-C., Celi, V., & Bulf, C. (2024). How to teach geometry in continuity along schooling? In T. Lowrie, A. Gutiérrez Rodrigues, & F. Emprin (Eds.), *Proceedings of the Twenty-Sixth ICMI Study. Advances in Geometry Education. 26th ICMI Study* (pp.231–238). IREM de Reims, 2024, 978- 2-910076-17-7. ⟨hal-04577863v2⟩

Salin, M.-H. (2002). Repères sur l’évolution du concept de milieu en théorie des situations [Key points on the evolution of the concept of milieu in the theory of situations]. In J.-L. Dorier, M. Artaud, M. Artigue, R. Berthelot, & R. Floris (Eds.), *Actes de la 11ème Ecole d’Eté de Didactique des Mathématiques* [Proceedings of the 11th summer school on the didactics of mathematics] (pp. 111–124). La Pensée Sauvage.